\documentclass[11pt]{amsart}
\usepackage{graphicx} 
\usepackage{setspace}

\usepackage[margin=1in]{geometry} 
\usepackage{amsmath,amssymb,amsfonts,amsthm}
\usepackage{mathrsfs}
\usepackage{array}
\usepackage{longtable}
\usepackage{graphicx}
\usepackage{float}
\usepackage{tikz}

\newtheorem{theorem}{Theorem}[section]

\newtheorem{corollary}[theorem]{Corollary}

\newtheorem{definition}{Definition}[section]
\newtheorem{example}{Example}
\newtheorem{remark}{Remark}

\newcommand{\rk}{\text{rk}}

\newcommand{\disc}{\text{disc}}
\newcommand{\NS}{\text{NS}}

\newcommand{\R}{\mathbb{R}}
\newcommand{\C}{\mathbb{C}}
\newcommand{\Q}{\mathbb{Q}}
\newcommand{\Z}{\mathbb{Z}}
\newcommand{\PP}{\mathbb{P}}

\newcommand{\bbm}{\begin{bmatrix}}
\newcommand{\ebm}{\end{bmatrix}}

\newcommand{\OO}{\mathcal{O}}

\newcommand{\Pic}{\text{Pic}}

\newcommand{\MWL}{\text{MWL}}
\newcommand{\Triv}{\text{Triv}}
\newcommand{\Nef}{\text{Nef}}

\newcommand{\NE}{\text{NE}}
\newcommand{\Amp}{\text{Amp}}

\makeatletter

\newcommand{\Rmnum}[1]{\expandafter\@slowromancap\romannumeral #1@}
\makeatother

\title{Newton-Okounkov polygons with a Small Number of Vertices and Picard number} 
\author{Yue Yu} 
\date{\today} 

\begin{document}

\begin{abstract}
        Newton--Okounkov bodies serve as a bridge between algebraic geometry and
convex geometry, enabling the application of combinatorial and geometric
methods to the study of linear systems on algebraic varieties.
This paper contributes to understanding the algebro-geometric information
encoded in the collection of all Newton--Okounkov bodies on a given
surface,
focusing on Newton--Okounkov polygons with few vertices and on elliptic K3
surfaces.

Let $S$ be an algebraic surface and $mv(S)$ be the maximum number of vertices of the Newton-Okounkov
bodies of $S$. We prove that $mv(S) = 4$ if and only if its Picard
number $\rho(S)$ is at least 2 and $S$ contains no negative irreducible
curve. Additionally, if $S$ contains a negative curve, then $\rho(S) =
2$ if and only if $mv(S) = 5$.

Furthermore, we provide an example involving two elliptic K3 surfaces
to demonstrate that when $\rho(S) \geq 3$, $mv(S)$ no longer determines
the Picard number $\rho(S)$.
        \end{abstract}

\maketitle

\section{Introduction}
       
        After their introduction by A. Okounkov in \cite{okounkov1996brunn},
        R.~Lazarsfeld and M.~Musta\c{t}\u{a} \cite{lazarsfeld2009convex} on one
        hand, and K.~Kaveh and A.~Khovanskii \cite{kaveh2012newton} on the other, independently
        developed the theory of Newton--Okounkov bodies of big divisors
        on algebraic varieties.
        Given a big divisor $D$ and a suitable flag $Y_\bullet$ on a projective
variety
        $X$, the Newton--Okounkov body $\Delta_{Y_\bullet}(D)$ collects the data
        of the orders of vanishing of sections of $\mathcal{O}(kD)$, $k\gg0$, along
        the components of $Y_\bullet$.
        In the case of surfaces, a rather precise description of
Newton-Okounkov bodies
        is known, in terms of Zariski decompositions,
        thanks to the seminal work by Lazarsfeld--Mustata, later developed
        by  A.~K\"uronya, V.~Lozovanu, and C.~Maclean \cite{kuronya2012convex},
        J.~Ro\'e and T.~Szemberg \cite{roe2022number} and others.
        All Newton--Okounkov bodies on surfaces are polygons whose area determines
        the \emph{volume} of $D$ (an asymptotic invariant measuring the growth
of $h^0(\mathcal{O}(kD))$, see [\cite{lazarsfeld2017positivity}, Section 2.2]).
        The \emph{shape} of the polygons is a more enigmatic invariant, but
        their number of vertices can be bounded in terms of the Picard
        number $\rho(X)$.
       
        J.~Ro\'e and T.~Szemberg in \cite{roe2022number} determined the maximal
        number of vertices of the Newton--Okounkov body of an ambple divisor,
$mv(S)$,
        in terms of the geometry of negative curves on an algebraic surface $S$. They proved that
        $\rho(S)=1$ if and only if $mv(S)=3$, that is, all
        Newton--Okounkov bodies of ample divisors are triangles.
        In this paper we extend their work characterizing surfaces with
$mv(S)\le 5$:
       
        \begin{theorem}
           Let \(S\) be an algebraic projective surface. Then:

    (1) \(mv(S) = 3\) if and only if \(\rho(S) = 1\).

    (2) \(mv(S) = 4\) if and only if \(\rho(S) > 1\) and there are no negative irreducible curves on \(S\).

    (3) If \(S\) contains some negative irreducible curve, then \(mv(S) = 5\) if and only if \(\rho(S) = 2\).  
        \end{theorem}

        Thus we see that when the surface has certain special properties,
        such as the presence of negative irreducible curves, $mv(S)$ can, to
some extent, characterize the Picard number if it is small.

        In the last part of the paper we initiate a study of Newton--Okounkov
bodies
        on elliptic K3 surfaces.
        These are very interesting surfaces in many respects, in particular for
        their configurations negative curves which will reflect in their
        Newton--Okounkov bodies.
        However, we show by an explicit example that even for elliptic K3 surfaces,
        if $\rho(S) > 2$, then $mv(S)$ does not determine the Picard number:
       
        \begin{theorem}
                There exist two elliptic K3 surfaces with different Picard numbers,
yet the maximal number of vertices of the Newton--Okounkov body of an
ample divisor on each surface is the same.
        \end{theorem}
       
        In conclusion, we show that the collection of Newton--Okounkov bodies
of ample
        divisors on a smooth projective surface contains information on the Picard
        number, beyond what was known in the literature. However, this data alone
        is insufficient to determine the Picard number,
        even under restrictive assumptions on the surface.
        Our approach for elliptic K3 surfaces suggests potential extensions
that could further link their geometry with the shape of their
Newton--Okounkov
        bodies.
        For more general classes of surfaces, we remark that in \cite{moyano2025newton} the authors conjecture that the Picard
        number of a surface should be determined by its Newton--Okounkov bodies
        \emph{if one allows for infinitesimal flags} (which are not considered
        in the present paper).
        Our results give partial evidence for this conjecture,  while also
confirming that, without infinitesimal flags, the conjecture does not hold.

        The paper is organized as follows. Section 2 introduces the preliminary
background for the proofs, including the method of counting the number
of vertices of Newton-Okounkov bodies proposed by J.~Ro\'e and T.~Szemberg,
        as well as fundamental aspects of lattice theory, elliptic K3 surfaces,
and the Mordell-Weil group. Section 3 discusses the cases where $mv(S) =
3, 4, 5$ and establishes their relationship with the Picard number.
Section 4 constructs an example to prove Theorem 4.1.

\subsection*{Conventions}
 This paper works over an algebraically closed field. All root lattices \(A_n\), \(D_n\), \(E_l\) are taken to be negative-definite. For an elliptic surface $f: S \to C$, we use $P, Q$, etc., to denote rational points in the generic fibre $E(k(C))$, while $(P), (Q)$, etc., represent their images as sections on the surface $S$, which are curves.

\section{Preliminaries}

\subsection{Newton-Okounkov bodies on surfaces}

Recall Okounkov's construction: Let \(X\) be a smooth irreducible projective variety of dimension \(d\). The construction depends on the choice of a fixed flag:
\[
Y_\bullet: X = Y_0 \supset Y_1 \supset \cdots \supset Y_d = \{\text{point}\},
\]
where \(Y_i\) is a smooth irreducible subvariety of codimension \(i\) in \(X\). Such flag is called admissible flag. Given a big divisor \(D\) on \(X\), one defines a valuation function:
\[
\nu = \nu_{Y_\bullet} = \nu_{Y_\bullet,D} : (H^0(X, \OO_X(D)) \setminus \{0\}) \to \Z^d, \quad s \mapsto \nu(s) = (\nu_1(s), \dots, \nu_d(s)).
\]
Set \(\nu_1 = \nu_1(s) = \operatorname{ord}_{Y_1}(s)\). Then \(s\) determines, in a natural way, a section \(\tilde{s}_1 \in H^0(X, \OO_X(D - \nu_1 Y_1))\) that does not vanish identically along \(Y_1\). By restricting, this gives a non-zero section \(s_1 \in H^0(Y_1, \OO_{Y_1}(D - \nu_1 Y_1))\). Then we take \(\nu_2(s) = \operatorname{ord}_{Y_2}(s_1)\), and continue in this manner to define the remaining \(\nu_i(s)\).

\noindent Next, define \(\nu(D) = \operatorname{Im}((H^0(X, \OO_X(D)) \setminus \{0\}) \to \Z^d)\) to be the set of valuation vectors of non-zero sections of \(\OO_X(D)\). It is not hard to check that \(\#\nu(D) = h^0(X, \OO_X(D))\). Define the Newton-Okounkov body of \(D\) with respect to the fixed flag \(Y_\bullet\) as:
\[
\Delta(D) = \Delta_{Y_\bullet}(D) = \text{closed convex hull}\left(\bigcup_{m \geq 1} \frac{1}{m} \cdot \nu(mD)\right).
\]
\(\Delta_{Y_\bullet}(D)\) is a convex body in \(\R^d\).

Let \(S\) be a smooth projective surface, and let \(D\) be any big divisor on \(S\). Consider an admissible flag \(Y_\bullet: S \supset C \supset \{p\}\). Consider the $\R$-divisor $D_t=D-tC$ and, if $D_t$ is effective or pseudo-effective, denote its Zariski decomposition
$$
D_t=P_t+N_t
$$

Let $\nu=\nu_C(D)$ be the coefficient of $C$ in the negative part $N_0$ of the Zariski decomposition of $D$, and let $\mu=\mu_C(D)=\max\{t\in\R\ |\ D_t\text{ is pseudo-effective}\}$. Note that $D_\mu$ belongs to the boundary of the pseudo-effective cone, it is not big. For every $t\in[\nu,\mu]$, define $\alpha(t)=(N_t\cdot C)_p$, that is the local intersection multiplicity of the negative part of $D_t$ and $C$ at $p$, and $\beta(t)=\alpha(t)+P_t\cdot C$. $\Delta_{Y_\bullet}(D)$ is the region in the plane $(t,s)$ defined by the inequalities $\nu\leq t\leq\mu$, $\alpha(t)\leq s\leq\beta(t)$ (\cite{lazarsfeld2009convex} Section b). Note that $\alpha$ and $\beta$ are continuous piecewise linear functions in the interval $[\nu,\mu]$, respectively convex and concave.

\noindent Paper \cite{roe2022number} showed an approach to understanding the number of vertices in Newton-Okounkov polygons by analyzing the dependence of \(N_t\) on \(t\), and from this derived information on the functions \(\alpha\) and \(\beta\). We directly reference several important definitions and Theorem \ref{thm01} from \cite{roe2022number} to compute the maximal number of vertices of the Newton-Okounkov body of an ample divisor on a surface.

\noindent Consider an effective divisor \(N = C_1 + \cdots + C_k\) which is a sum of irreducible curves on a projective surface \(S\), whose intersection matrix is negative definite. According to \cite{roe2022number}, one defines two numbers associated with \(N\): Let \(mc(N)\) denote the largest number of irreducible components of a connected divisor contained in \(N\). And let
    \begin{equation}
        mv(N) = 
        \begin{cases} 
            k + mc(N) + 4, & \text{if } k < \rho - 1, \\
            k + mc(N) + 3, & \text{if } k = \rho - 1.
        \end{cases}
    \end{equation}

\noindent On surface $S$, if we take all \(N\), one can define:
\[
mv(S) = \max\{mv(N) \mid N = C_1 + \cdots + C_k \text{ negative definite}\}.
\]

\noindent There is the following theorem from\cite{roe2022number}:

\begin{theorem}[\cite{roe2022number} Theorem 5.5]\label{thm01}
    Let \(S\) be any smooth projective surface, and let \(D\) be a big divisor on \(S\). Then,
    \begin{equation}
        \max_{Y_\bullet} \{\# \text{vertices}(\Delta_{Y_\bullet}(D))\} \leq mv(S),
    \end{equation}
    where the maximum is taken over all admissible flags \(Y_\bullet\). Moreover, if \(D\) is ample, then for every \(3 \leq \nu \leq mv(S)\), there exists a flag \(Y_\bullet\) such that \(\Delta_{Y_\bullet}(D)\) has exactly \(\nu\) vertices.
\end{theorem}

This implies that maximum \(mv(S)\) is always attainable. And the number \(mv(S)\) is well-defined and bounded above by \(2\rho + 1\) by \cite{kuronya2012convex}.
\subsection{Lattices}

In this section, we introduce the fundamental concept of lattices. The main reference is Ebeling \cite{ebeling2013lattices}.

A \textbf{lattice} \(L\) is a finitely generated free \(\Z\)-module of finite rank endowed with a non-degenerate symmetric bilinear form \((\cdot, \cdot)\) with values in \(\Z\). A lattice \(L\) is called \textbf{even} if \((x, x) \in 2\Z\) for all \(x \in L\); otherwise, \(L\) is called \textbf{odd}. The determinant of the intersection matrix with respect to an arbitrary basis over \(\Z\) is called the \textbf{discriminant}, \(\disc(L)\).

A lattice \(L\) and the \(\R\)-linear extension of its bilinear form \((\cdot, \cdot)\) give rise to the real vector space \(L_{\R} := L \otimes_\Z \R\), endowed with a symmetric bilinear form. \(L_\R\) can be diagonalized with only \(1\) and \(-1\) on the diagonal. The \textbf{signature} of \(L\) is \((n_+, n_-)\), where \(n_+\) (respectively \(n_-\)) is the number of \(+1\) (resp. \(-1\)) on the diagonal. The lattice \(L\) is called \textbf{definite} if either \(n_+ = 0\) or \(n_- = 0\); otherwise, \(L\) is \textbf{indefinite}. 

\begin{example}
The hyperbolic plane is the lattice 
\begin{equation}
    U := 
    \begin{bmatrix}
        0 & 1 \\
        1 & 0
    \end{bmatrix}.
\end{equation}
That is, \(U \cong \Z^2 = \Z \cdot e \oplus \Z \cdot f\) with the quadratic form given by \((e^2) = (f^2) = 0\) and \((e \cdot f) = 1\). The discriminant \(\disc U = -1\).

The hyperbolic lattice has significant applications in elliptic K3 fibrations. The intersection matrix formed by the zero section and general fiber is 
\begin{equation}
    \begin{bmatrix}
        -2 & 1 \\
        1 & 0
    \end{bmatrix} \cong U .
\end{equation}
The determinant of this matrix is \(-1\), and its signature is \((1, 1)\). It is isometric to the hyperbolic lattice.
\end{example}

To a lattice \(L\), we can define the \textbf{dual lattice} as \(L^* = \{x \in L \otimes \Q \mid (x, L) \in \Z\}\). \(L^*\) is a lattice with the extended bilinear product on \(L\), and it has the same rank as \(L\). Regarding the discriminant of \(L\), we have \(\disc(L) = |L^*/L|\). A lattice \(L \subset \R^n\) is called \textbf{unimodular} if \(L^* = L\).

Let \(\Gamma\) be a lattice in \(\R^n\). A \(\Z\)-submodule \(\Lambda\) of \(\Gamma\) is called a sublattice of \(\Gamma\). It is a lattice in some subspace \(W \subset \R^n\), which is isomorphic to \(\R^k\) for some \(k\). A sublattice \(\Lambda\) of \(\Gamma\) is called \textbf{primitive} if \(\Gamma / \Lambda\) is a free \(\Z\)-module.

Let \(L\) be an even, negative-definite lattice. Define \(R := \{x \in \Gamma \mid x^2 = -2\}\). An element \(x \in R\) is called a \textbf{root}. \(L\) is called a \textbf{root lattice} if \(R\) generates \(\Gamma\).

Every root lattice is the orthogonal direct sum of irreducible root lattices associated with Coxeter-Dynkin diagrams. For example, the Coxeter-Dynkin diagram corresponding to \(A_n\) is showed below (representing each \(e_i\) by a vertex, and connecting them by edges if \(e_i \cdot e_j = 1\)):

\begin{figure}[H]
    \centering
    \begin{tikzpicture}
        \fill (0, 0) circle (2pt);
        \fill (1, 0) circle (2pt);
        \fill (2, 0) circle (2pt);
        \fill (3, 0) circle (2pt);
        \fill (4, 0) circle (2pt);
        \fill (5, 0) circle (2pt);
        \draw (0, 0) -- (1, 0);
        \draw (1, 0) -- (2, 0);
        \draw[dotted] (2, 0) -- (3, 0);
        \draw (3, 0) -- (4, 0);
        \draw (4, 0) -- (5, 0);
    \end{tikzpicture}
    \caption{\(A_n\) (\(n\) vertices)}
\end{figure}
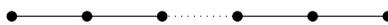

\subsection{K3 Surfaces}

A $\textbf{K3 surface}$ over $\C$ is a complete non-singular surface $X$ such that $\omega_X\cong\OO_X$ and $H^1(X,\OO_X)=0$. Some general references are \cite{huybrechts2016lectures} \cite{barth2015compact}.

$H^2(X,\Z)$ is endowed with a unimodular intersection pairing, making it isometric to the K3 lattice
$$
\Lambda=U^3\oplus E^2_8,
$$

\noindent where $U$ is the hyperbolic lattice and $E_8$ is the unique up to isometry even unimodular negative definite lattice of rank $8$. The signature of $H^2(X,\Z)$ is $(3,19)$.

We denote by $\Pic (X)$ the group of invertible sheaves up to isomorphism, which is isomorphic to the group of divisors modulo linear equivalence. The $\textbf{Néron-Severi group}$ $\text{NS}(X)$ is the group of divisors modulo algebraic equivalence, it is a sublattice of $H^2(X,\Z)$, has signature $(1,\rho(X)-1)$, here $\rho(X)=\text{rk}(\text{NS}(X))$ is the $\textbf{Picard number}$ of $X$.

\begin{example}
    A smooth quartic $X\subset\PP^3$ is a K3 surface.
\end{example}

\subsection{Elliptic Surfaces}

Assume $C$ is a smooth projective curve over $k$. Let $S$ be a smooth projective algebraic surface and
$$
f:S\to C
$$

\noindent is a surjective morphism with connected fibres defined over $k$. Then $f:S\to C$ is called a $\textbf{fibration}$ on the surface $S$ with the base curve $C$. For each point $v\in C(k)$, let $F_v=f^{-1}(v)$ be the $\textbf{fibre}$ (as a divisor) over $v$. A fibration is $\textbf{relatively minimal}$ means that no fibre contains a $(-1)$-curve. We suppose that $S$ is relatively minimal. If all fibres $F_v$ except for finitely many ones are smooth curves of genus $g$ we call $f:S\to C$ a $\textbf{genus $g$ fibration}$. The $\textbf{generic fibre}$ of $f$ is the fibre over the generic point of the base curve $C$. Thus the generic fibre is a curve of genus $g$ defined over the function field $k(C)$ of $C$ over $k$. A $\textbf{section}$ of $f$ is a morphism $\sigma:C\to S$ such that $f\circ \sigma:C\to C$ is the identity map of $C$. For each fibre $F_v(v\in C)$, $\sigma(v)$ is a point on $F_v$.

\begin{definition}\label{def311}
    An $\textbf{elliptic surface}$ is a genus $1$ fibration
    $$
    f:S\to C
    $$
    from a smooth projective surface $S$ to a smooth projective curve $C$ with a section $\sigma_0:C\to S$ which is relatively minimal.
\end{definition}

For an elliptic surface, we can choose the section $\sigma_0:C\to S$ as the $\textbf{zero section}$ such that all but finitely many fibres $F_v$ are elliptic curves with the origin $\sigma_0(v)$. 

Assuming that $F_v=f^{-1}(v)$ is a singular fiber where $v\in C(k)$, we can express it as a divisor on $S$ with multiplicities as follows:
$$
F_v=\sum_{i=0}^{m_v-1}\mu_{v,i}\Theta_{v,i}
$$

\noindent here $m_v$ is the number of the irreducible components in $F_v$; $\Theta_{v,i}\ (0\leq i\leq m_v-1)$ are irreducible components; $\mu_{v,i}$ is the multiplicity of $\Theta_{v,i}$ in $F_v$.

$\textbf{Zariski's lemma}$ states that if $F_v$ is a reducible singular fibre $(m_v>1)$, then every component $\Theta_{v,i}$ of $F_v$ is a smooth rational curve which has self-intersection number $(\Theta_{v,i})^2=-2$.

According to Kodaira's classification \cite{kodaira1963compact}, all possible types of reducible singular fibers are classified into the following types with $m>1$ and $b\geq 0$:
\begin{align}\label{sf}
    \Rmnum{1}_m,\ \Rmnum{1}^*_b,\ \Rmnum{3},\ \Rmnum{4},\ \Rmnum{2}^*,\ \Rmnum{3}^*,\ \Rmnum{4}^*.
\end{align}
\noindent In this paper, we only consider singular fibres of type $\Rmnum{1}_2$ and $\Rmnum{3}$. They can be described as follows:  
\begin{itemize}
    \item A singular fibre of type $\Rmnum{1}_2$ consists of two smooth rational curves intersecting at two points with multiplicity $1$.
    \item A singular fibre of type $\Rmnum{3}$ consists of two smooth rational curves intersecting at a single point with multiplicity $2$.
\end{itemize}

For detailed representations of these fibres, readers may refer to Figures 5.1 and 5.2 in \cite{schutt2019mordell}, where the structure of the aforementioned singular fibres is illustrated.

\subsection{Mordell-Weil Lattice}

The main reference for this section is Chapter 6 of \cite{schutt2019mordell}. For additional details, readers may consult \cite{miranda1989basic}, \cite{shioda1990mordell}, and the foundational work \cite{cox1979intersection}.

The \textbf{trivial lattice} $\Triv(S)$ (\cite{schutt2019mordell}, Definition 6.1) is the sublattice of $\NS(S)$ generated by the zero section and fibre components:
\[
\Triv(S) = \langle (O), F \rangle \bigoplus_{\nu \in R} T_\nu^-,
\]
where $R$ denotes the finite subset $\{\nu \in C(k) \mid T_\nu \neq 0\}$ of points on the base curve where singular fibres are located. The rank of the trivial lattice is given by \cite{schutt2019mordell}, Proposition 6.3:
\[
\rk(\Triv(S)) = 2 + \sum_{\nu \in R} (m_\nu - 1).
\]

For every elliptic surface $S$, the Néron-Severi lattice $\NS(S)$ is finitely generated and torsion-free. Moreover, algebraic and numerical equivalence on $S$ coincide (\cite{schutt2019mordell}, Theorem 6.4).

Let $f: S \to C$ be an elliptic surface with generic fibre $E$ over $K = k(C)$. For every rational point on the generic fibre, the map $P \mapsto (P) \ \mod{\Triv(S)}$ defines an isomorphism of abelian groups $E(K) \cong \NS(S)/\Triv(S)$ (\cite{schutt2019mordell}, Theorem 6.5). If $S$ has a singular fibre, then the Mordell-Weil group $E(K)$ is a finitely generated abelian group (\cite{schutt2019mordell}, Theorem 6.6). Let $r = \rk(E(K))$. The following equality holds (\cite{schutt2019mordell}, Corollary 6.7):
\begin{equation}\label{shioda formula}
    \rho(S) = r + 2 + \sum_{\nu \in R} (m_\nu - 1).
\end{equation}

For any $P \in E(K)$, there exists a unique element in $\NS(S)_{\mathbb{Q}}$, denoted by $\varphi(P)$. $\varphi(P)$ is the class of the divisor $D_P$ with $\mathbb{Q}$-coefficients. When $S$ has no reducible singular fibre, $D_P = (P) - (O) - ((P) \cdot (O) + \chi)F$, where $(P), (O)$ are sections and $F$ is the divisor of a fibre.

For any $P, Q \in E(K)$, the height pairing $\langle P, Q \rangle$ is defined as:
\[
\langle P, Q \rangle = \chi + (P) \cdot (Q) + (Q) \cdot (O) - (P) \cdot (Q) - \sum_{\nu \in R} \text{contr}_\nu(P, Q),
\]
and, in particular,
\[
\langle P, P \rangle = 2\chi + 2(P) \cdot (O) - \sum_{\nu \in R} \text{contr}_\nu(P).
\]
Here, $\text{contr}_\nu(P, Q)$ represents the local contribution from the singular reducible fibre at $\nu \in R$. In this discussion, $\text{contr}_\nu(P, Q)$ is always zero. Readers interested in a complete definition can consult \cite{schutt2019mordell}, Lemma 6.16.

Thus, this defines a $\mathbb{Q}$-valued symmetric bilinear pairing on $E(K)$, inducing a positive-definite lattice structure on $E(K) / E(K)_{\text{tors}}$. The pairing on the Mordell-Weil group $E(K)$ is called the \textbf{height pairing}, and the lattice $(E(K)/E(K)_{\text{tors}}, \langle \cdot, \cdot \rangle)$ is referred to as the \textbf{Mordell-Weil lattice (MWL)} of the elliptic surface $f: S \to C$.

In the case where the elliptic surface $S$ has no torsion section and only irreducible singular fibres, its Néron-Severi lattice can be expressed as:
\begin{equation}\label{nsequ}
    \NS(S) = \Triv(S) \oplus \MWL(S)^-.
\end{equation}

\section{The case of $mv(S)=3,4,5$}

\begin{theorem}\label{thm1}

    Let \(S\) be an algebraic projective surface. Then:

    (1) \(mv(S) = 3\) if and only if \(\rho(S) = 1\).

    (2) \(mv(S) = 4\) if and only if \(\rho(S) > 1\) and there are no negative irreducible curves on \(S\).

    (3) If \(S\) contains some negative irreducible curve, then \(mv(S) = 5\) if and only if \(\rho(S) = 2\).
    
    \begin{proof}

    As explained in the introduction, the case \(mv(S) = 3\) is well known. For an algebraic surface \(S\), \(mv(S) = 3\) if and only if \(\rho(S) = 1\). 

    Proof of (2): Suppose there exists \(N = C_1 + \dots + C_k\) such that \(mv(S) = mv(N) = k + mc(N) + s = 4\). This implies the only possible case is that \(s = 4\) and \(k = 0\). Then \(\rho(S) > k + 1 = 1\). Conversely, if \(\rho(S) > 1\) and there are no negative irreducible curves on \(S\), then \(N = \emptyset\). Thus, \(mv(S) = 4\).  

    Proof of (3): \textbf{“Only if”}: Suppose \(S\) has Picard number \(2\). From the assumption, it contains a negative curve \(C\). Let \(N = C\). Then \(k = mv(N) = 1\) satisfies \(k = \rho - 1\), \(mv(N) = k + mc(N) + 3 = 5\). As \(mv(N) = 2\rho + 1 = 5\), this is the maximum possible value of \(mv(S)\). Therefore, \(mv(S) = mv(N) = 5\), as desired.

    \textbf{“If”}: Suppose \(mv(S) = 5\). Then there exists \(N = C_1 + \dots + C_k\) such that \(mv(S) = k + mc(N) + s = 5\), where \(s = 3\) if \(\rho(S) = k + 1\), and \(s = 4\) if \(\rho(S) > k + 1\). Clearly, \(k \neq 0\); otherwise, \(mv(S) < 5\). As \(k \geq 1\), it follows that \(1 \leq mc(S) \leq k\). Thus, \(k + mc(N) \geq 2\). This implies \(s = 3\), and it follows that \(k = mc(N) = 1\). In this case, \(\rho = k + 1 = 2\).

    \end{proof}
\end{theorem}

\begin{corollary}\label{cor1}
    Let \(S\) be an elliptic surface with a section such that its Euler characteristic \(\chi(S) \neq 0\). Then the Picard number \(\rho(S) = 2\) if and only if the maximum number of vertices of the Newton-Okounkov bodies of an ample divisor \(mv(S)\) is \(5\).

    \begin{proof}
        By the assumption \(\chi(S) \neq 0\), and it can only be positive. Thus, the zero section \((O)\) has self-intersection number \((O)^2 = -\chi(S) < 0\) \cite{schutt2019mordell}. By Theorem \ref{thm1}, the result follows.
    \end{proof}
\end{corollary}

\begin{remark}
    The condition \(\chi(S) \neq 0\) in the above corollary is equivalent to the absence of singular fibers. It is worth noting that the condition \(\chi(S) \neq 0\) cannot be omitted. Consider the trivial elliptic surface \(S = \PP^1 \times E\), where \(E\) is an elliptic curve. In this case, since there are no singular fibers, \(\chi(S) = 0\). The effective cone \(\NE(S)\) is generated by \(E\) and \((O) \cong \PP^1\) (\cite{lazarsfeld2017positivity} Example 1.4.33), both of which have self-intersection \(0\). In this case, \(N = \emptyset\), and \(mv(S) = 4\).
\end{remark}

\begin{corollary}
    Let \(S\) be an elliptic surface with a section, and suppose \(mv(S) = 4\). Then its Euler characteristic \(\chi(S) = 0\).

    \begin{proof}
        It is obvious that \(\rho(S) \geq 2\). Since \(mv(S) = 4\), this implies that there are no negative irreducible curves on this surface. This requires the self-intersection number of the section \((O)^2=-\chi(S) \geq 0\), and there are no singular fibers. Hence, the Euler characteristic \(\chi(S) = 0\).
    \end{proof}
\end{corollary}

\begin{example}
    Consider the abelian surface \(A = E \times E\), where \(E\) is an elliptic curve without complex multiplication. The Picard number of \(A\) is \(3\), $\overline{NE}(A)=\Nef(A)$ (\cite{lazarsfeld2017positivity} Lemma 1.5.4), so $mv(A)=4$ and there are no irreducible negative curves. Furthermore, \(\chi(S) = 0\).
\end{example}

\section{Counterexample}

Let $S$ be a K3 surface. For $\rho(X) \geq 3$, either $S$ does not contain any smooth rational curves at all, or $\overline{\NE}(S)$ is the closure of the cone spanned by all smooth rational curves $C \subset S$ (\cite{huybrechts2016lectures}, Chapter 8, Section 3). In the case of elliptic K3 surfaces, there is always a smooth rational curve, namely the zero section. Hence, the latter case holds, $\overline{\NE}(S) \neq \overline{\Amp}(S)$, and $\overline{\NE}(S)$ has no circular parts.

Therefore, we can restrict our attention to $-2$-curves and the intersection lattice with respect to $N=C_1+\dots+C_k$ can be expressed as the orthogonal sum of finitely many root lattices.

\begin{theorem}\label{thm4}

    There exist two elliptic K3 surfaces with different Picard numbers $\rho_1,\rho_2>2$, yet the maximal number of vertices of the Newton-Okounkov body of an ample divisor on each surface is the same.

\begin{proof}
    Consider the first K3 surface \(S_1\) with Picard number \(\rho(S_1) = 3\) such that it admits an elliptic fibration with exactly one reducible singular fiber of type \(\mathrm{I}_2\). According to the Picard number formula \eqref{shioda formula}: 
    \[
    \rho = 3 = r + 2 + \sum_{\nu \in R} (m_\nu - 1),
    \]
    we have \(r = 0\). Additionally, we assume there is no torsion section, so by equation \ref{nsequ}, the Néron-Severi group of \(S_1\) is \(\mathrm{NS}(S_1) = U \oplus A_1\). By (\cite{huybrechts2016lectures} Chapter 14 Corollary 3.1), such \(S_1\) with Néron-Severi group $\NS(S_1)$ exists. 

    Next, we compute \(mv(S_1)\). Let \(N_1\) consist of the zero section and one of the components of the reducible singular fiber \(\mathrm{I}_2\) that intersects the zero section at a single point with multiplicity $1$. The intersection matrix of \(N_1\) is \(A_2\), which is negative definite. Therefore,
    \[
    mv(S_1) =k+mc(N_1)+3= 2 + 2 + 3 = 7.
    \]
    This is also the maximum value achievable for a surface with Picard number \(3\), as \(2\rho + 1 = 7\).

    Next, consider another K3 surface \(S_2\) with Picard number \(\rho(S_2) = 4\) such that \(S_2\) admits an elliptic fibration with no reducible singular fibers. From the formula \eqref{shioda formula}:
    \[
    \rho = 4 = r + 2 + \sum_{\nu \in R} (m_\nu - 1) = r + 2,
    \]
    we deduce \(r = 2\). Assume these two nonzero sections, denoted by \((P)\) and \((Q)\), do not intersect each other or the zero section. We compute the Mordell-Weil lattice (\(MWL\)) as follows. 

    Let \(P\) and \(Q\) represent the rational points on elliptic curve over function field. The pairing in \(MWL\) is computed as:
    \begin{equation}
        \langle P, Q \rangle = \chi + (P)\cdot (O) + (Q)\cdot (O) - (P)\cdot(Q) = \chi = 2,
    \end{equation}
    \begin{equation}
        \langle P, P \rangle = 2\chi + 2(P)\cdot (O) = 2\chi = 4 = \langle Q, Q \rangle.
    \end{equation}
    By equation \ref{nsequ}, the Néron-Severi group of \(S_2\) is:
    \[
    \mathrm{NS}(S_2) = U \oplus 
    \begin{bmatrix}
    -4 & -2 \\
    -2 & -4
    \end{bmatrix}.
    \]
    The signature of \(\mathrm{NS}(S_2)\) is \((1, 3)\). By (\cite{huybrechts2016lectures} Chapter 14 Corollary 3.1), such a K3 surface with this Néron-Severi lattice exists. Let the basis of \(\mathrm{NS}(S_2)\) be the classes of \((O), F, D_P, D_Q\) (using the same symbols where there is no ambiguity), where \((O), F\) are the basis of \(U\) with intersection relations \((O)^2 = -2\), \((O) \cdot F = 1\), and \(F^2 = 0\).

We now compute \(mv(S_2)\) using lattice methods. First, we construct an embedding \(3A_1 \hookrightarrow \mathrm{NS}(S_2)\). Let the basis of \(3A_1\) be \((O), (P), (Q)\), and \((P), (Q)\) can be expressed in terms of the basis of \(\mathrm{NS}(S_2)\) as:
\[
(P) = D_P + (O) + 2F, \quad (Q) = D_Q + (O) + 2F.
\]

    \[
    (P) = D_P + (O) + 2F, \quad (Q) = D_Q + (O) + 2F.
    \]

    Thus, \(N_2\) can be constructed as three disjoint negative curves. Therefore,
    \[
    mv(N_2) = k+mc(N_2)+3= 3 + 1 + 3 = 7.
    \]

    Finally, we prove that no \(N_3\) exists such that \(mv(N_3) \geq 8\). If \(mv(N_3) = 8\), the possible configurations are \(mv(N_3) = 3 + 2 + 3 = 8\) or \(mv(N_3) = 2 + 2 + 4 = 8\). If \(mv(N_3) = 9\), the only configuration is \(mv(N_3) = 3 + 3 + 3 = 9\). 

In all cases, there must exist at least two negative curves intersecting each other. Since \(S_2\) is a K3 surface, the negative irreducible curves must be \((-2)\)-curves (\cite{bauer2004zariski} section 3.2). To ensure \(N_3\) is negative definite, the two \((-2)\)-curves must intersect at exactly one point. Suppose such two curves exist. Then, the lattice \(A_2\) must be embedded into the intersection lattice of \(N_3\). Consequently, \(A_2 \hookrightarrow \mathrm{NS}(S_2)\).

    Let the basis of \(A_2\) be \(E_1, E_2\), and
    \[
    E_1 = aD_P + bD_Q + c(O) + dF, \quad E_2 = a'D_P + b'D_Q + c'(O) + d'F.
    \]
    Here \(a, b, c, d, a', b', c', d'\) are integers satisfying:
    \[
    (aD_P + bD_Q + c(O) + dF) \cdot (a'D_P + b'D_Q + c'(O) + d'F) = 1,
    \]
    \[
    (aD_P + bD_Q + c(O) + dF)^2 = (a'D_P + b'D_Q + c'(O) + d'F)^2 = -2.
    \]
    Expanding these equations yields the system:
    \begin{align}
        2a^2 + 2b^2 + 2ab + c^2 - cd &= 1, \label{eq1} \\
        2a'^2 + 2b'^2 + 2a'b' + c'^2 - c'd' &= 1, \label{eq2} \\
        -2(2aa' + ab' + a'b + 2bb' + cc') + cd' + dc' &= 1. \label{eq3}
    \end{align}
    Next verify if this system admits integer solutions. 

From equation \eqref{eq1}, the term \(2a^2 + 2b^2 + 2ab\) is even. Hence, for \(c^2 - cd = c(c-d)\), we have \(2 \nmid c(c-d)\), implying that $c(c-d)$ is odd. Thus, \(2 \nmid c\) and \(2 \nmid (c-d)\), which further implies \(2 \mid d\).

From equation \eqref{eq3}, the term \(-2(2aa' + ab' + a'b + 2bb' + cc')\) is even. Therefore, \(2 \nmid (cd' + dc')\). Without loss of generality, suppose \(2 \mid cd'\) and \(2 \nmid dc'\). Then, \(2 \nmid d\neq0\) and \(2 \nmid c'\neq0\). However, according to the earlier conclusion, \(2 \mid d\), leading to a contradiction. Similarly, if we assume \(2 \nmid cd'\) and \(2 \mid dc'\), we also arrive at a contradiction.

Therefore, \(A_2 \not\hookrightarrow \mathrm{NS}(S_2)\). This proves that \(mv(S_2) = 7 = mv(S_1)\), yet the two surfaces have different Picard numbers.

\end{proof}
\end{theorem}

\bibliographystyle{plainurl}
\bibliography{NOB}

\end{document}